\documentclass[a4paper]{article}
\usepackage{amsmath} % because of {split}. amssymb is not needed now?
\usepackage{amsfonts}
\usepackage{latexsym}
\newcommand{\ep}{\hspace*{\fill}$\Box$}
\newcommand{\eps}{\varepsilon}
\newcommand{\pr}{{\bf Proof. }}

\newcommand{\cl}{\mbox{\rm cl}}

\newcommand{\R}{\mathbb R}
\newcommand{\N}{\mathbb N}
\newcommand{\C}{\mathbb C}

\newcommand{\rmd}{\mbox{\rm d}}
\newcommand{\amtil}{\ensuremath{\tilde{\mathcal{A}}_m(X)} } %changed from M to X, mike 27/02/01
\newcommand{\io}{\iota}

\newcommand{\gd}{\ensuremath{{\mathcal G}^d} }

\newcommand{\comp}{\subset\subset}

\newcommand{\cinfty}{{\cal C}^\infty}

\newtheorem{thr}{\hspace*{-3mm} \bf}[section]
\newcommand{\bt}{\begin{thr} {\bf Theorem. }}
\newcommand{\et}{\end{thr}}
\newcommand{\bp}{\begin{thr} {\bf Proposition. }}
\newcommand{\bc}{\begin{thr} {\bf Corollary. }}
\newcommand{\blem}{\begin{thr} {\bf Lemma. }}
\newcommand{\bex}{\begin{thr} {\bf Example. }\rm} 
\newcommand{\bexs}{\begin{thr} {\bf Examples. }\rm}
\newcommand{\bd}{\begin{thr} {\bf Definition. }}
\newcommand{\beast}{\begin{eqnarray*}}
\newcommand{\eeast}{\end{eqnarray*}}

\newcommand{\rem}[1]{\vadjust{\rlap{\kern\hsize\thinspace\vbox%
                       to0pt{\hbox{${}_\clubsuit${\small\tt #1}}\vss}}}}
\newcommand{\ahat}{\ensuremath{\hat{\mathcal{A}}_0(X)} }
\newcommand{\atil}{\ensuremath{\tilde{\mathcal{A}}_0(X)} }
 
\newcommand{\ehat}{\ensuremath{\hat{\mathcal{E}}(X)} } 
\newcommand{\emhat}{\ensuremath{\hat{\mathcal{E}}_m(X)} }
\newcommand{\nhat}{\ensuremath{\hat{\mathcal{N}}(X)} }
\newcommand{\ghat}{\ensuremath{\hat{\mathcal{G}}(X)} } 
\newcommand{\lhat}{\ensuremath{\hat{L}_\zeta} }                    
\newcommand{\al}{\alpha}
\newcommand{\bet}{\beta} 

\newcommand{\Om}{\Omega}\newcommand{\om}{\omega}
\newcommand{\si}{\sigma}\newcommand{\la}{\lambda}
\newcommand{\vphi}{\varphi}

\newcommand{\D}{{\cal D}}

\newcommand{\pa}{\partial}
 \newcommand{\G}{{\cal G}}

\newcommand{\supp}{\mathop{\mathrm{supp}}}

\newcommand{\A}{{\mathcal A}}

\newcommand{\be}{ \begin{equation} }\newcommand{\ee}{\end{equation} }
\newcommand{\beq}{ \begin{equation} }\newcommand{\eeq}{\end{equation} }
\newcommand{\bea}{\begin{eqnarray}}\newcommand{\eea}{\end{eqnarray}}
\newcommand{\beas}{\begin{eqnarray*}}\newcommand{\eeas}{\end{eqnarray*}}
\newcommand{\beqs}{\begin{equation*}}\newcommand{\eeqs}{\end{equation*}}
\newcommand{\lb}{\label}

\newcommand{\ben}{\begin{enumerate}}\newcommand{\een}{\end{enumerate}}
\newcommand{\ba}{\begin{array}}\newcommand{\ea}{\end{array}}
\newcommand{\bthm}{\begin{thr} {\bf Theorem. }}
\newcommand{\bprop}{\begin{thr} {\bf Proposition. }}
\newcommand{\bcor}{\begin{thr} {\bf Corollary. }}
\newcommand{\bdef}{\begin{thr} {\bf Definition. }}
\newcommand{\brem}{\begin{thr} {\bf Remark. }\rm}
\newcommand{\bth}{\begin{thr}\rm}
\newcommand{\ethi}{\end{thr}}

\newcommand{\lgl}{\langle}
\newcommand{\rgl}{\rangle}
\newcommand{\ca}{{\cal A}}

\newcommand{\cc}{{\cal C}}
\newcommand{\cd}{{\cal D}}

\title{Diffeomorphism invariant Colombeau algebras.\newline Part III: Global theory}
\author
{ Michael Kunzinger\footnote{Universit\"at Wien, Institut f\"ur Mathematik, Strudlhofgasse
4, Austria;  electronic mail: {\tt Michael.Kunzinger@univie.ac.at}}}
\date{%\today; preliminary version
}

\begin{document}
\maketitle

\begingroup

\begin{abstract}
We present the construction of an associative, commutative algebra $\ghat$
of generalized functions on a manifold $X$ satisfying the following optimal 
set of permanence properties:
\begin{itemize}
\item[(i)] $\D'(X)$ is linearly embedded into $\ghat$, $f(p)\equiv 1$
is the unity in $\ghat$.
\item[(ii)] For every smooth vector field $\xi$ on $X$ there exists a derivation
operator $\hat L_\xi: \ghat \to \ghat$ which is linear and satisfies
the Leibniz rule.
\item[(iii)] $L_\xi|_{\D'(X)}$ is the usual Lie derivative.
\item[(iv)] $\circ |_{\cinfty(X)\times\cinfty(X)}$ is the pointwise product
of functions.
\end{itemize}
Moreover, the basic building blocks of $\ghat$ are defined in purely 
intrinsic terms of the manifold $X$.
\end{abstract}
\section{Introduction}
Recent applications of Colombeau's theory of algebras of generalized functions
to problems of a primarily geometric nature (\cite{clarke}, \cite{geo}, \cite{geo2}, 
\cite{penrose}, \cite{vickersESI}, \cite{symm}, \cite{DKP}) have very clearly indicated
the need for a global intrinsic version of the construction on differentiable
manifolds. The development of local diffeomorphism invariant Colombeau algebras 
on open subsets of $\R^n$, initiated in \cite{CM, Jel}, has only recently 
been completed
in \cite{found}. Based on \cite{found, vim} as well as Parts I and II
of this series (also in this volume) in the present article we present
an intrinsic global version of Colombeau's theory on differentiable manifolds
adapted to the needs of applications in mathematical physics. 

In what follows we shall use freely notation and terminology from \cite{found}.
Additionally, we will use the following conventions. $X$ will denote
an orientable, paracompact, $n$-dimensional smooth manifold.
with atlas $\mathfrak{A} = \{(U_\alpha,\psi_\alpha)\ :\
\alpha \in A\}$. $\Omega^n_c(X)$ is the space of compactly
supported (smooth) $n$-forms on $X$. For coordinates $y^1,\dots,y^n$
on $U_\al$, elements of $\Om_c^n(\psi_\al(U_\al))$ will be written as
$\vphi\, d^ny := \vphi\, dy^1 \wedge \dots \wedge dy^n$.  
% The pullback of any $\om\in \Omega^n_c(\psi_\al(U_\al))$ under $\psi_\al$ is
% denoted by $\psi_\al^*(\om)$.  Then $\int_X \psi_\al^* (\vphi\,d^ny) =
% \int_{\psi_\al(U_\al)} \vphi\,d^ny$ for all $\vphi\,d^ny \in
% \Omega^n_c(\psi_\al(U_\al))$. 
% We generally
% use the following convention: 
If $B\subseteq \psi_\al(U_\al)$ then
$\hat B := \psi_\al^{-1}(B)$ and if $f: \psi_\al(U_\al)\to \R$ resp.\ $\C$ then
$\hat f := f\circ\psi_\al$.  
% A similar ``hat-convention''
% will be applied to the function spaces to be defined in the following
% section. Moreover, since $X$ is supposed to be oriented we can and shall
% identify $n$-forms and densities henceforth.
% For an open subset $\Om$ of $\R^n$, the space of distributions on
% $\Om$ (i.e., the dual of the (LF)-space $\D(\Om)$) will be denoted by
% $D'(\Om)$.  
% For a diffeomorphism $\mu: \tilde\Om \to \Om$, the
% pullback of any $u\in \D'(\Om)$ under $\mu$ is defined by
% \begin{equation} \label{distpull}
% \langle \mu^*(u),\varphi\rangle = \langle u(y),
% \varphi(\mu^{-1}(y))\cdot|\mathrm{det}\, D\mu^{-1}(y)|\rangle
% \end{equation}
Concerning distributions on manifolds we follow the terminology
of~\cite{deR} and \cite{Marsden}. The space of distributions on
$X$ is defined as $\D'(X) = \Omega^n_c(X)'$ (the dual of the space
of compactly supported $n$-forms) and  
% Observe that in this
% setting test objects no longer have {\em function} character but are
% $n$-forms.  
% In the context of Colombeau algebras it is natural to regard smooth
% functions as regular distributions (which in the non-orientable case
% enforces the use of test {\em densities}; this is in accordance 
% with~\cite{Hoe} but has to be distinguished from the setting 
% of~\cite{D3}). 
operations on distributions are defined as (sequentially) continuous extensions of
classical operations on smooth functions.
For example, for $\zeta\in \mathfrak{X}(X)$ (the space of smooth vector fields 
on $X$) and $u\in \D'(X)$ the
Lie derivative of $u$ with respect to $\zeta$ is given by $\langle L_\zeta u,
\om\rangle = -\langle u, L_\zeta \om \rangle$. If $u\in \D'(X)$,
$(U_\al,\psi_\al) \in \mathfrak{A}$ then the local representation of
$u$ on $U_\al$ is the element $(\psi_\al^{-1})^*(u)\in
\D'(\psi_\al(U_\al))$ defined by
\begin{equation}  \label{locdist}
\langle (\psi_\al^{-1})^*(u), \vphi \rangle = \langle u|_{U_\al},
\psi_\al^*(\vphi d^n y) \rangle \qquad \forall \vphi \in
\D(\psi_\al(U_\al))
\end{equation}
\section{Testing procedures} %moderateness, negligibility
An intrinsic formulation on differentiable manifolds of the diffeomorphism 
invariant Colombeau algebra $\gd$ introduced in \cite{found}
faces a number of serious obstacles due to the following indispensible 
technical ingredients of the construction of $\gd$: 
\begin{itemize}
\item Translation (convolution) used for the embedding of $\D'$
into $\gd$, leading to terms of the form 
$\langle u, \vphi(. - x)\rangle$
($u\in \D'(\Om)$, $\vphi\in\D(\Om)$, $\Om\subseteq \R^n$ open).
\item Scaling operations of the form $\vphi \to S_\eps\vphi
:= \frac{1}{\eps}\vphi(\frac{.}{\eps})$.
\item Moment integrals of the form $\int \xi^\beta \vphi(\xi)\, d\xi$.
\end{itemize}
Obviously, none of these operations allows a direct generalization to the
manifold setting. Our task therefore consists in unfolding the
diffeomorphism invariant `essence' underlying the local 
testing procedures used for determining relationship in the spaces
of moderate resp.\ negligible mappings (cf.\ \cite{found}, ch.\ 3).
Definition \ref{kernels} below  introduces test objects on $X$
(so called {\em smoothing kernels}) 
that display precisely those properties of local test objects 
corresponding to regularization via convolution and linear scaling on $\R^n$
(part (i)) resp.\ to the interplay between $x$- and $y$-differentiation
in the local context.

\bd
\beas
 \hat{{\mathcal A}}_0(X) &:=& \{\om \in \Omega^n_c(X)\ : \ \int \om = 1 \} \\
% \]
% \et
% Starting with the essential of {\bf (D1)},
% we choose the basic space for the forthcoming definition of the Colombeau algebra on $X$
% as follows:
% \bd
 \hat{{\mathcal E}}(X) &:=& \cc^\infty(\hat{{\mathcal A}}_0(X) \times
X)
\eeas
\et
$\hat{{\mathcal E}}(X)$ is our basic space, both moderate and negligible
maps will be elements of $\hat{{\mathcal E}}(X)$. Localization of elements
of $R\in\hat{{\mathcal E}}(X)$ is effected by the map
\[
(\psi_\al^{-1})^* R := R \circ  (\psi_\al^* \times \psi_\al^{-1}) 
\in \hat{\mathcal E}(\psi_\al(U_\al))
\]
% is  an  element of $\hat{\mathcal E}(\psi_\al(U_\al))$. More generally, if $\mu:
% X_1  \to  X_2$ is a diffeomorphism and $R\in \hat\E(X_2)$ then its pullback
% $\hat\mu(R)  \in  \hat\E(X_1)$  under  $\mu$  is defined as $R\circ \bar \mu$
% where $\bar\mu(\om,p) = ((\mu^{-1})^*\om,\mu(p))$ and clearly $(\mu_1 \circ
% \mu_2)\hat{} = \hat\mu_2 \circ \hat\mu_1$ for diffeomorphisms $\mu_1$, $\mu_2$.
Suppose that $f:X\times X \to \bigwedge^n T^*X$ is smooth such that 
% for each $(p,q)\in X\times X$,
% $f(p,q)$ belongs to the fiber over $q$ or, equivalently, that 
for each $p\in X$,
$f_p := (q\mapsto f(p,q))$ is contained in $\Om^n(X)$. 
% Obviously, $p\mapsto f_p  
% \in \cc^\infty(X,\Om^n(X))$ in this case. 
Then for any $\zeta \in {\mathfrak X}(X)$, we introduce the following
two notions of Lie derivatives of $f$ with respect to $\zeta$.
% which, essentially, 
% arise as Lie derivatives of $q \mapsto f(p,q)$ resp.\ $p\mapsto f(p,q)$: 
% On the one hand, viewing $\Om^n(X)$ together with the topology of pointwise (i.e., 
% fiberwise) convergence as a locally convex space we define
\beas
(L'_\zeta f)(p,q)&:=& L_\zeta(p \mapsto f(p,q)) = \left. \frac{d}{dt} 
\right|_0f((\mathrm{Fl}^\zeta_t)(p),q)\\
(L_\zeta f)(p,q) &:=& L_\zeta(q \mapsto f(p,q)) = \left. \frac{d}{dt} 
\right|_0((\mathrm{Fl}^\zeta_t)^*f_p)(q)
\eeas
After these preparations we finally turn to the definition of smoothing kernels,
the global analogue on $X$ of the translated scaled test objects 
$T_xS_\eps\phi(\eps,x)$ in the local theory: 

\bd \label{kernels}
$\Phi\in \cc^\infty(I \times X,\ahat) %\subseteq \cc^\infty(I\times X\times X,\Lambda^nT^*X)
$  
is called a smoothing kernel if %it satisfies the following conditions
\begin{itemize}
  \item [(i)] $\forall K\subset\subset X$ $\exists\, \eps_0$, $C>0$ $\forall p\in K$
              $\forall \eps \le \eps_0$:  $\supp\Phi(\eps,p)\subseteq B_{\eps C}(p)$
  \item [(ii)] $\forall K \subset\subset X$ $\forall k,\, l\in \N_0$ 
        $\forall \zeta_1,\dots,\zeta_k,\theta_1,\dots,\theta_l\in\mathfrak{X}(X)$ 
        \[
                  \sup_{{p\in K}\atop{q\in X}} \|L_{\theta_1}\dots
                  L_{\theta_l}(L'_{\zeta_1}+L_{\zeta_1})\dots
            (L'_{\zeta_k}+L_{\zeta_k})
                        \Phi(\eps,p)(q) \| = O(\eps^{-(n+l)})
                \]
\end{itemize}
The space of smoothing kernels on $X$ is denoted by \atil.
\et
Here $B_{\eps C}(p)$ denotes the ball of radius $\eps C$, measured with
respect to the distance function $d_h$ induced on $X$ by some 
Riemannian metric $h$ on $X$ and $\|\,.\,\|$ denotes the norm induced 
by $h$ on $\bigwedge^n T^*X$. Both notions are independent of the chosen
metric $h$ (cf.\ \cite{vim}, Lemma 3.4).

The grading on ${\mathcal A}_0(\R^n)$ into the subspaces ${\mathcal A}_m(\R^n)$
consisting of those test functions of unit integral whose moments up to order
$m$ vanish is replaced by the following subspaces of $\tilde{\mathcal A}_0(X)$
in the global case.
\bd \label{smoothingkernels}
Let $\tilde{\mathcal A}_m(X)$ ($m\in \N$) be
the set of all $\Phi\in \atil$ satisfying
\[
\forall f\in \cc^\infty(X) \ \forall K\subset\subset X \
\sup_{p\in K}  |f(p)- \int_X f(q)\Phi(\eps,p)(q) | = O(\eps^{m+1})
\]
\et
The technical motivation for the exact form of this definition is that it provides
precisely what is needed for proving $\iota\big|_{\cc^\infty} = \sigma$ later on.
Moreover, in the local case the above requirement in essence reproduces the original 
spaces ${\mathcal A}_m$.
It is shown in \cite{vim}, Lemma 3.7 that the spaces $\tilde{\mathcal A}_m(X)$
are in fact nontrivial for all $m\in\N_0$.
% Next, we introduce the appropriate notion of Lie derivative for elements of 
% $\hat{\mathcal E}(X)$.
\bd For any $R\in \ehat$ and any $\zeta\in \mathfrak{X}(X)$ we define
the Lie derivative of $R$ with respect to $\zeta$ by
\begin{equation} \label{liealg}
(\hat{L}_\zeta R)(\om,p) := - \rmd_1R(\om,p)(L_\zeta \om) + L_\zeta(R(\om,\,.\,))|_p
\end{equation}
where $\rmd_1 R(\om,x)$ denotes the derivative of $\om \to R(\om,x)$ (cf.\
\cite{KM}).
\et
In order to derive the local form of this Lie derivative we 
let $\zeta\in \mathfrak{X}(X)$, set $\zeta_\alpha = (\psi_\al^{-1})^*(\zeta|_{U_\al})$
% Since the derivative of the (linear and continuous) map $\psi_\al^*: 
% \Om_c^n(\psi_\al(U_\al)) \to
% \Om_c^n(U_\al)$ in any point equals 
% the map itself, on $U_\al$ we obtain (writing $R$ in place of $R|_{U_\al}$): 
and calculate as follows
\begin{eqnarray*}
\lefteqn{
((\psi_\al^{-1})^*(\hat L_\zeta R))(\vphi\,d^ny,x)}\\
&&\hphantom{}
=(\hat L_\zeta R)(\psi_\al^*(\vphi\,d^ny),\psi_\al^{-1}(x))\\
&&\hphantom{}
= L_\zeta(R(\psi_\al^*(\vphi\,d^ny),\,.\,))(\psi_\al^{-1}(x))\\ 
&&\hphantom{mmmmmmmm}
- (\rmd_1R)(\psi_\al^*(\vphi\,d^ny),\psi_\al^{-1}(x))
(\underbrace{L_\zeta(\psi_\al^*(\vphi\,d^ny))}_{\psi_\al^*(L_{\zeta_\alpha}(\vphi\,d^ny))}) \\
&&\hphantom{}
= [L_{\zeta_\alpha}((\psi_\al^{-1})^* R))](\vphi\,d^ny,x) - 
[\rmd_1((\psi_\al^{-1})^* R)(\vphi\,d^ny,x)](L_{\zeta_\alpha}(\vphi\,d^ny))
\end{eqnarray*}
Now setting $\zeta_\al=\pa_{y^i}$ ($1\le i \le n$) we obtain the local algebra 
derivative $D_i^J$ in the J-formalism (cf.\ \cite{found}, Ch.\ 5 and \cite{Jel}). 
%\et
\section{The algebra $\ghat$}
We begin by introducing the subspaces of moderate and negligible maps
of $\hat {\cal E}(X)$. 
% Finally we are in a position to define the subspaces of 
% moderate and negligible 
% elements of $\ehat$.  
\bd     \lb{defmodghat}
%\mbox{\rm{\bf (D3)}}\ \ $R\in \ehat$ is moderate if
% $\forall K\subset\subset X$  $\forall k\in \N_0$  $\exists N\in \N$
%  $\forall\ \zeta_1,\dots,\zeta_k\in
%                \mathfrak{X}(X)$ 
%  $\forall\ \Phi \in \atil$ 
%\begin{equation} \label{mod}
%\sup_{p\in K} |L_{\zeta_1}\dots L_{\zeta_k}(R(\Phi(\eps,p),p)) | = O(\eps^{-N})
%\end{equation}
$R\in \ehat$ is called {\rm moderate} if
the following condition is satisfied:
$$\forall K\subset\subset X\  \forall k\in \N_0\ \exists N\in \N\
\forall\ \zeta_1,\dots,\zeta_k\in\mathfrak{X}(X)\ 
\forall\ \Phi \in \atil$$
\begin{equation} \label{mod}
\sup_{p\in K} |L_{\zeta_1}\dots L_{\zeta_k}(R(\Phi(\eps,p),p)) | =
O(\eps^{-N}).
\end{equation}

The subset of moderate elements of \ehat is denoted by \emhat.
\et
\bd   \lb{defnegghat}
%\mbox{\rm{\bf (D4)}}\ \ $R\in \ehat$ is called negligible if it satisfies
%\\
%%\parbox{11cm}{
%\begin{equation}\label{*}
%\begin{array}{l}
%\forall K\subset\subset X \ \forall k, l \in \N_0 \ \exists m \in \N
%    \ \forall \  \zeta_1,\dots,\zeta_k\in
%                \mathfrak{X}(X)    \\
%\forall \Phi \in  \tilde{\mathcal A}_m(X) \ 
%       \sup_{p\in K} |L_{\zeta_1}\dots L_{\zeta_k}(R(\Phi(\eps,p),p)) | 
%       = O(\eps^{l})          
%\end{array}
%\end{equation}
%%} \hfill \parbox{8mm}{\bea \label{*}\eea}\\ 
$R\in \ehat$ is called {\rm negligible} if
the following condition is satisfied:
%\parbox{11cm}{
$$\forall K\subset\subset X\  \forall k,l\in \N_0\ \exists m\in \N\
\forall\ \zeta_1,\dots,\zeta_k\in\mathfrak{X}(X)\ 
\forall\ \Phi \in \amtil$$
\begin{equation} \label{*}
\sup_{p\in K} |L_{\zeta_1}\dots L_{\zeta_k}(R(\Phi(\eps,p),p)) | =
O(\eps^{l}).
\end{equation}
The subset of negligible elements of \ehat will be denoted by \nhat.
\et
We shall see in \ref{idealcharghat} below that in fact Lie derivatives can
be omitted completely in the definition of \nhat if we additionally suppose
that $R\in \emhat$. This fact is another instance of a very general result 
(Th.\ 13.1 of \cite{found}) stating a similar characterization of the 
Colombeau ideal as a subspace of the space of moderate functions without 
resorting to derivatives for practically all types of Colombeau algebras.
We immediately obtain
%The following result is immediate from the definitions:
\bt
\begin{itemize}
\item[(i)]  \emhat is a subalgebra of \ehat. 
\item[(ii)]  \nhat is an ideal in \emhat.
\end{itemize}
%\ep
\et
In order to connect the development of $\ghat$ already at this early stage
to that of $\gd$ we need a localization procedure allowing to transport smoothing
kernels on the manifold to local test objects on open subsets of $\R^n$ and
vice versa. However, a direct translation is not feasible since localizations
of smoothing kernels display rather poor properties concerning domain of 
definition. To precisely formulate the translation process we introduce the
following spaces of functions. We denote by $C_b^\infty(I\times\Om,\A_0(\R^n))$
the space of smooth maps $\phi: I\times \Omega \to \A_0(\R^n)$ 
such that for each $K \comp \Om$
(i.e., $K$ a compact subset of $\Om$) and any $\al\in\N_0^n$, the set 
$\{\pa^\al_x\phi(\eps,x)\mid\eps\in(0,1],x\in  K\}$
is bounded in $\cd(\R^n)$. For any $m\ge 1$ we set
\beas 
\A_m^\Box(\Om) 
&:=& 
\{\Phi \in C_b^\infty(I\times\Om,\A_0(\R^n)) 
\mid|\sup_{x\in K}|\int \Phi(\eps,x)(\xi)\xi^\al d\xi| =  \\
&& \hphantom{\A_m^\Box(\Om) := \ mmmmmg\,} O(\eps^{m})  \ 
(1 \le |\al| \le m) \ \forall
K\comp \Om \} \\[6pt] 
\A_m^\Delta(\Om) &:=& \{\Phi \in
C_b^\infty(I\times\Om,\A_0(\R^n)) \mid|\sup_{x\in K}|\int
\Phi(\eps,x)(\xi) \xi^\al d\xi| =  \\
&& \hphantom{\A_m^\Box(\Om) := \ mmm}
% \{\Phi \in
% C_b^\infty(I\times\Om,\A_0(\R^n)) |} 
O(\eps^{m+1-|\al|}) \ (1 \le |\al| \le m) \ \forall
K\comp \Om \}.
\eeas
Elements of $\A_m^\Box(\Om)$ are said to have
asymptotically vanishing moments of order $m$ (more precisely, in the
terminology \cite{found} elements of $\A_m^\Box(\Om)$ are of type
$[\mathrm{A}_g]$, the abbreviation standing for {\em asymptotic} vanishing of
moments {\em globally}, i.e., on each $K\comp \Om$).  
% The spaces  $\A_m^\Box(\Om)$
% and $ \A_m^\Delta(\Om) $ 
% %play a crucial role in the characterizations 
% %of    the    algebra   (\cite{found},   Section 10) and in
% %Section~\ref{localization} below (cf.\  also the discussion of diffeomorphism
% %invariance at the end of this subsection).
% will play a crucial role in the proof of localization of $\ghat$ to
% $\gd(\psi_\al(U_\al))$; this, in turn, will be an indispensable
% ingredient for deriving {\bf (T1)}, {\bf (T4)}, {\bf (T5)}, thus
% for the construction of $\ghat$ as a differential algebra
% containing $\cd'(X)$ as a whole.
It is easily seen that $\A_m^\Box(\Om)\subseteq \A_m^\Delta(\Om) \subseteq
\A_{2m-1}^\Box(\Om)$.

Moreover, by 
%$\cc^\infty_b(I\times\Om,\ca_0(\R^n))$
%of test objects used for constructing $\gd(\Om)$ 
%we define
$\cc^\infty_{b,w}(I\times\Om,\ca_0(\R^n))$ we denote the space of all
$\phi: D \to
\ca_0(\R^n)$  where $D$ is some subset (depending on $\phi$) of
$(0,1]\times\Om$ such that for $D,\phi$ the following holds: 

For each
$L\subset\subset\Om$ there exists $\eps_0$ and a  subset $U$
of $D$ which is open in $(0,1]\times\Om$ such that
\bea 
&& \hspace*{-3em}
(0,\eps_0]\times L\subseteq U(\subseteq D) \ \mbox{and} \ \phi \ 
\mbox{is smooth on } \ U \label{atoz1}\\
&& \hspace*{-3em}
\{\pa^\bet\phi(\eps,x)\mid0<\eps\le\eps_0,\ x\in L\} \ \mbox{is bounded in} \
\cd(\R^n) \  \forall \,\bet\in\N_0^n \label{atoz2}
\eea
Here the  subscript  $w$ indicates the  weaker  requirements  on  the  domain 
of definition  of $\phi$.
% Note that $\cc^\infty_{b,w}(I\times\Om,\ca_0(\R^n))$ precisely
% consists
% of those $\phi$ satisfying the requirements listed in condition
% $(\mathrm{Z})$ of \rf{a--z}.
The subspace of $\cc^\infty_{b,w}(I\times\Om,\ca_0(\R^n))$ consisting of those
$\phi$ whose moments up to order $m$ vanish asymptotically on each compact
subset of $\Om$ will be written as $\A_{m,w}^\Box(\Om)$.
$\A_{m,w}^\Delta(\Om)$ is defined analogously. 
With this terminology at hand we can now precisely state the transport properties
of smoothing kernels.

\blem\label{local_kernels} 
Denote by $(U_\al,\psi_\al)$ a chart in $X$. 

{\em (A) Transforming smoothing kernels to local test objects.}
\begin{itemize}
\item[(i)] Let $\Phi$ be a smoothing kernel. Then the map $\phi$ 
defined by
\begin{equation} \label{lockeri}
\phi(\eps,x)(y)d^ny:=\eps^n\,((\psi_\al^{-1})^*\Phi(\eps,\psi_\al^{-1}x))
    (\eps y+x) 
\end{equation}
$(x\in\psi_\al(U_\al),\,y\in\R^n)$ is an element of 
$\cc^\infty_{b,w}(I\times\psi_\al(U_\al),\ca_0(\R^n))$.

\item[(ii)] If, in addition, $\Phi \in \tilde{{\mathcal A}}_{m}(X)$ for some 
$m\in \N$ then $\phi\in{{\mathcal A}}^{\Delta}_{m,w}(\psi_\al(U_\al))$, i.e., 
\begin{equation} \label{n2}
\int \phi(\eps,x)(y)y^\bet dy = O(\eps^{m+1-|\bet|})
\quad (1\le |\bet| \le m)
\end{equation}
uniformly on compact sets. In particular, if $\Phi \in 
\tilde{{\mathcal A}}_{2m-1}(X)$  then $\phi \in
{\mathcal A}_{m,w}^{\Box}(\psi_\al(U_\al))$.
\end{itemize}

{\em (B) Transporting local test objects onto the manifold.}
\begin{itemize}
\item[(i)] Let $\phi \in \cc_b^\infty(I\times\Om,\ca_0(\R^n))$ 
and $\Phi_1\in\atil$. Let $K\subset\subset\psi_\al(U_\al)$, 
$\chi,\chi_1\in\D(\psi_\al(U_\al))$ with $\chi\equiv 1$ on an open neighborhood 
of $K$ and $\chi_1\equiv 1$ on an open neighborhood of
$\supp\chi$. Then
\begin{equation}\label{lockerii}
\begin{array}{rcl}
\Phi(\eps,p) \!\!\!&:=&\!\!\!
    (1-\hat\chi(p)\la(\eps))\Phi_1(\eps,p) \\[0.3em]
    &&+\hat\chi(p)\la(\eps)\psi_\al^*\left(\frac{1}{\eps^n}\phi(\eps,\psi_\al p)
    (\frac{y-\psi_\al p}{\eps})\chi_1(y)d^ny\right) 
\end{array}
\end{equation}
%
%   das ist versuch vom 01 03 20    m.g.         
%
%\begin{eqnarray}
%\lefteqn{\Phi(\eps,p):=
%    (1-\hat\chi(p)\la(\eps))\Phi_1(\eps,p)}\nn\\
%    &&\hphantom{mm}+\hat\chi(p)\la(\eps)\psi_\al^*\left(\frac{1}{\eps^n}\phi(\eps,\psi_\al p)
%    \left(\frac{y-\psi_\al p}{\eps}\right)
%    \chi_1(y)d^ny\right)\quad\quad  \label{lockerii} 
%\end{eqnarray}
%
%   ende versuch
%
%where $\hat\chi$ denotes  $\chi\circ\psi_\al$
is a smoothing kernel (with $\la$ a smooth cut-off function). 
\item[(ii)] If, in addition,
$\phi\in{{\mathcal A}}^{\Delta}_{m}(\psi_\al(U_\al))$
and $\Phi_1\in \amtil$ then $\Phi\in \amtil$. 
In particular, if 
$\phi \in {\mathcal A}_m^{\Box}(\psi_\al(U_\al))$ and $\Phi_1\in \amtil$ then
$\Phi\in \amtil$.
\end{itemize}
\et

For a proof of this result we refer to \cite{vim}.

Using the above result we may now derive local characterizations of moderateness
and negligibility, thereby establishing the link to $\gd$ promised above.
\bt \label{modloc}
%{\rm(Localization of moderateness)} 
Let $R\in \ehat$. Then
for all $\al$, 
\[
R\in \emhat \Leftrightarrow (\psi_\al^{-1})^* (R|_{U_\al}) \in
{\mathcal E}_M(\psi_\al(U_\al)).
\]
\et
\pr
($\Rightarrow$) Let $R\in \emhat$, $(U_\al,\psi_\al) \in {\mathfrak A}$ 
and set $R' := (\psi_\al^{-1})^* (R|_{U_\al})$. 
Then for $K\comp \psi_\al(U_\al)=:V_\al$, $\bet \in \N_0^n$ 
and $\phi\in$ $\cc^\infty_b(I\times\psi_\al(U_\al),$ $\ca_0(\R^n))$. 
we have to show that
\begin{equation}\label{modloc1}
\sup_{x\in K}|\pa^\bet (R'(T_xS_\eps\phi(\eps,x),x))| = 
O(\eps^{-N})
\end{equation} 
for some $N\in \N$. To this end we fix some 
$\Phi_1\in\atil$ and define $\Phi\in \atil$ by (\ref{lockerii}). 
Let $\pa^\beta = \pa_{i_1}\dots \pa_{i_{|\bet|}}$ and for
$1\le i_j \le n$ choose $\zeta_{i_j}\in \mathfrak{X}(X)$ such that the local expression
of $\zeta_{i_j}$ coincides with $\pa_{i_j}$ on a neighborhood of $K$. Then for $\eps$ 
sufficiently small (\ref{modloc1}) equals
\[
\sup_{p\in  \hat K}|L_{\zeta_{i_1}}\dots L_{\zeta_{i_{|\bet|}}} R(\Phi(\eps,p),p)| \, .
\]

($\Leftarrow$) Let  $\zeta_1,\dots,\zeta_k\in \mathfrak{X}(X)$ and 
%(without loss of generality)
$\hat K \comp U_\al$ for some chart $(U_\al,\psi_\al)$. Let $\Phi\in \atil$ and
define $\phi$ by (\ref{lockeri}). Since $\phi: D \to \ca_0(\R^n)$ belongs to
$\cc^\infty_{b,w}(I\times\psi_\al(U_\al),\ca_0(\R^n))$ it follows from
\cite{found}, Th.\ 10.5
that given $\beta\in \N_0^n$ there exists some $N\in
\N$ and some $\eps_0>0$ such that for $x\in K$ and $\eps\le \eps_0$ we have
$(\eps,x)\in D$ and $|\pa^\beta (\psi_\al^{-1})^*
(R|_{U_\al})(T_xS_\eps\phi(\eps,x),x))| = 
O(\eps^{-N})$. From this and (\ref{mod}) the result follows.
% Inserting this into the local representation of (\ref{mod})
% immediately gives the result. 
\ep

\bt \label{locneg}\lb{negloc} %\hspace*{-5pt}{\rm (Localization of negligibility)} 
Let $R \in \emhat$. \nolinebreak 
\hspace*{-5pt} Then for all $\al$,
\[
R\in \nhat \Leftrightarrow (\psi_\al^{-1})^* (R|_{U_\al}) \in
\mathcal{N}(\psi_\al(U_\al)) .
\]
\et
\pr
($\Rightarrow$) Let $R\in \emhat$, $(U_\al,\psi_\al)\in {\mathfrak A}$ 
%some chart in $X$ 
and set
$R' := (\psi_\al^{-1})^* (R|_{U_\al})$. 
Let $K\comp \psi_\al(U_\al)=:V_\al$, $l\in \N$ and 
$\bet \in \N_0^n$.
Set $k=|\bet|$ and choose $m\in \N_0$  such that 
\begin{equation} \label{locneg1} 
   \sup_{p\in \hat K} |L_{\zeta_1}\dots L_{\zeta_k}(R(\Phi(\eps,x),x)) | = O(\eps^{l})          \end{equation}
for all $\zeta_1,\dots,\zeta_k \in \mathfrak{X}(X)$ and all $\Phi \in
\tilde{\mathcal A}_m(X)$. Then by \ref{local_kernels}, for 
$\phi\in {\mathcal A}_m^{\Box}(V_\al)$,  
$\Phi$ as in (\ref{lockerii}) is an element of $\Phi\in \tilde{\mathcal A}_m(X)$. 
% By
% \ref{local_kernels}, $\Phi\in \tilde{\mathcal A}_m(X)$, 
Hence (with $\zeta_{i_j}$ as in the proof of \ref{modloc}) we obtain
%$$
%\begin{array}{rcl}
%\sup_{x\in K}|\pa^\bet (R'(T_xS_\eps\phi(\eps,x),x))| &=& \\[0.2em]
%\sup_{p\in \hat K} |L_{\zeta_{i_1}}\dots L_{\zeta_{i_{|\bet|}}}(R(\Phi(\eps,x),x))|
%&=& O(\eps^l)
%\end{array}
%$$
\beas
\lefteqn{\sup\limits_{x\in K}|\pa^\bet (R'(T_xS_\eps\phi(\eps,x),x))|}\\
&&\hphantom{mmmmmmm}=\sup\limits_{p\in \hat K}
|L_{\zeta_{i_1}}\dots L_{\zeta_{i_{|\bet|}}}(R(\Phi(\eps,x),x))|
= O(\eps^l)
\eeas
 
so the claim follows from the characterization results in \cite{found} (Ch.\ 7 and 10).
% the characterization of negligibility following the
% definition of $\cn(\Om)$ (section \ref{notterm3}).

($\Leftarrow$) Let $k\in \N_0$, $l\in \N$, $\zeta_1,\dots,\zeta_k\in \mathfrak{X}(X)$
and $\hat K \comp U_\al$. %By \ref{negrem}
%the discussion at the end of section \ref{notterm3} 
There exists $m'$ such that %(with $R' = (\psi_\al^{-1})^* R|_{U_\al}$)
\begin{equation} \label{locneg2}
\sup_{x\in K}|\pa^\bet (R'(T_xS_\eps\phi(\eps,x),x))| = 
O(\eps^l)
\end{equation}
for all $|\bet|\le k$ and all $\phi \in {\mathcal A}_{m',w}^{\Box}(V_\al)$. Now set
$m = 2m'-1$ and let $\Phi \in \tilde {\mathcal A}_m(X)$. Then $\phi$ defined
by (\ref{lockeri}) is in ${\mathcal A}_{m',w}^{\Box}(V_\al)$ by
\ref{local_kernels} (A) (ii).
Hence inserting local representations of $\zeta_1,\dots,\zeta_k$, (\ref{locneg2}) 
implies (\ref{locneg1}).
\ep\medskip\\

From this and \cite{found}, Th.\ 13.1 we conclude
\bc \lb{idealcharghat}
 Let $R\in \hat{\mathcal E}_M(X)$. Then
\[
R\in \nhat \Leftrightarrow (\ref{*}) \ \mathrm{ holds \ for } \ k=0.
\]
\et
%\pr This follows directly from \ref{locneg} by taking into account
%\cite{found}, Th.\ 7.13 and \cite{found}, Th.\ 13.1. 
%\rf{JT1823}.\ep\medskip\\
As a final ingredient in the construction of $\ghat$, in the following
result we establish 
stability of \emhat and \nhat under Lie derivatives 
\bt \label{liestab} Let $\zeta\in \mathfrak{X}(X)$. Then
\begin{itemize}
\item[(i)] $\lhat \emhat \subseteq \emhat$.
\item[(ii)]  $\lhat \nhat \subseteq \nhat$.
\end{itemize}
\et
\pr
Let $R\in \emhat$, $\zeta\in \mathfrak{X}(X)$. By \ref{modloc} for any chart
$(U_\al,\psi_\al)$ we have $(\psi_\al^{-1})^* 
(R|_{U_\al}) \in {\mathcal E}_M(\psi_\al(U_\al))$.
Thus by the stability of the space of moderate functions on $\R^n$
(\cite{found}, Th.\ 7.10)
also $L_{\zeta_\al}(\psi_\al^{-1})^* (R|_{U_\al})$ $=$
$(\psi_\al^{-1})^* (\lhat R|_{U_\al}) \in {\mathcal E}_M(\psi_\al(U_\al))$
(where $\zeta_\al$ denotes the local representation of $\zeta$), which, 
again by \ref{modloc} gives the result. 
The claim for \nhat follows analogously from the stability of the space of
negligible functions under differentiation (\cite{found}, Th.\ 7.11).
\ep\medskip\\
%Finally we are in a position to define our main object of interest:
Having collected all the necessary properties of $\emhat$ and $\nhat$
(cf.\ the general scheme of construction for Colombeau algebras given in
\cite{found}, Ch.\ 3) we turn to the actual definition of $\ghat$.
\bd The full Colombeau algebra on $X$ is defined as
\[
\ghat := \emhat \big/ \nhat
\]
\et
$\ghat$ is a differential algebra with respect to the Lie derivative $\hat L$
induced by (\ref{liealg}), cf.\ \ref{liestab}.
For any $R \in \ehat$, its class in $\hat\G(X)$ will
% as usual 
be denoted by $\cl[R]$.

% %\newpage
% \lb{embedghat}
\section{Embedding properties}

In order to complete the list of properties of $\ghat$ given in the
abstract of this paper we still need to embed $\cinfty(X)$ and $\D'(X)$
into $\ghat$. 
\bd
Let  $u\in \D'(X)$, $f\in\cinfty(X)$;
\beas
(\io u)(\om,x)&:=&\lgl u,\om\rgl\\
(\si f)(\om.x)&:=&f(x).
\eeas
\et
$\io:\cd'(X)\to\ehat$ and $\si:\cinfty(X)\to\ehat$ are linear
embeddings, $\si$ respecting multiplication and unit
of $\cinfty(X)$. 

Next we show that $\io$ and $\si$ commute with
arbitrary Lie derivatives. Concerning $\io$, we obtain
from (\ref{liealg})
\beas
\hat L_\zeta(\iota u)(\om,p)
     &=& -\rmd_1 (\io u)(\om,p)(L_\zeta\om)+ 
%        \underbrace{
        L_\zeta((\io u)(\om,\,.\,))|_p
%       }_{=0}
                                  \\
     &=&-\rmd_1\left((\om,p)\mapsto\lgl
                 u,\om\rgl\right)(L_\zeta\om)\ +\ 0\\
     &=&-\lgl u,L_\zeta\om\rgl\\
     &=&\iota(L_\zeta u)(\om,p).
\eeas
Similarly, for $\si$ the first term  $-\rmd_1 (\si f)(\om,p)(L_\zeta\om)$
of $\hat L_\zeta(\si f)(\om,p)$ vanishes, so
$L_\zeta((\si f)(\om,\,.\,))|_p=L_\zeta f(p)=\si(L_\zeta f)(\om,p)$.

The most important properties of $\sigma$ and $\iota$ are collected in the
following result.
\bt
\begin{itemize}
\item[(i)] $\iota(\D'(X)) \subseteq \emhat$. 
\item[(ii)] $\sigma(\cinfty(X)) \subseteq \emhat$. 
\item[(iii)] $(\iota-\sigma)(\cinfty(X)) \subseteq \nhat$. 
\item[(iv)] $\iota(\D'(X))\cap \nhat = \{0\}$. 
\end{itemize}
\et
\pr
% {\bf (T1)}, finally, essentially follows from the local result \rf{th1} via
% \rf{modloc} resp.\ \rf{negloc}. Part (ii)
% ($\si(\cinfty(X))\subseteq\emhat$) being obvious, let us 
To show (i),
%i.e., $\io(\cd'(X))\subseteq\emhat$. To this end,
let $\om \in \D(\psi_\al(U_\al))$; then 
\begin{eqnarray*}
&&((\psi_\al^{-1})^* ((\io u)|_{U_\al}))(\om,x) = ((\io u)|_{U_\al})(\psi_\al^*(\om\,
d^ny),\psi_\al^{-1}(x)) = \\
&& \langle u, \psi_\al^*(\om\, d^ny)\rangle = 
\langle (\psi_\al^{-1})^*(u|_{U_\al}),\om\rangle
\end{eqnarray*} 
Since $(\psi_\al^{-1})^*(u|_{U_\al}) \in \D'(\psi_\al(U_\al))$ it follows from
the local theory that 
indeed $(\psi_\al^{-1})^* ((\io u)|_{U_\al}) \in {\mathcal E}_M(\psi_\al(U_\al))$.

(ii) is immediate and (iii) follows from \ref{locneg} and the 
corresponding local result.

Finally, to establish (iv) suppose
that $\io u \in \nhat$. By the same reasoning as above  
$(\om,x) \to \langle(\psi_\al^{-1})^*(u|_{U_\al}),\om\rangle 
\in \mathcal{N}(\psi_\al(U_\al))$ 
for each $\al$. Thus again by the respective local result
$(\psi_\al^{-1})^*(u|_{U_\al})=0$ for each $\al$, i.e. $u=0$.
\ep%\ms
%Summing up and passing to the quotient $\ghat$, we have 
\bc\label{lastcol}
$\iota: \D'(X) \to \ghat$ is a linear embedding that commutes
with Lie derivatives and coincides with $\sigma: \cc^\infty(X)\to \ghat$ 
on $\cc^\infty(X)$. Thus $\iota$ renders $\D'(X)$ a linear
subspace  and $\cc^\infty(X)$ a faithful subalgebra of \ghat.  
\et

\ref{lastcol} completes the construction of an intrinsic Colombeau algebra on
$X$ preserving all the distinguishing properties of the local algebra $\gd$.

%\newpage
% \bibliographystyle{plain}
% \bibliography{diana}

\end{document}